\documentclass[11pt,a4paper,twoside,draft]{amsart}
\usepackage{amsmath,bbm,amssymb,amsxtra}
\usepackage{enumerate}
\allowdisplaybreaks

\theoremstyle{definition}
\newcommand{\sptext}[3]{\hspace{#1 em}\mbox{#2}\hspace{#3 em}}
\newtheorem{lemma}{Lemma}[section] 
\newtheorem{proposition}[lemma]{Proposition}
 
\newtheorem{theorem}[lemma]{Theorem} 
\newtheorem{corollary}[lemma]{Corollary}
\newtheorem{remark}[lemma]{Remark}

\newtheorem{assumption}[lemma]{Assumption} 

\newcommand{\prop}[1]{\begin{proposition}\label{#1}\sl }
\newcommand{\eprop}{\end{proposition}}
\newcommand{\thm}[1]{\begin{theorem}\label{#1}\sl }
\newcommand{\ethm}{\end{theorem}}
\newcommand{\lem}[1]{\begin{lemma}\label{#1}\sl }
\newcommand{\elem}{\end{lemma}}

\newcommand{\beqno}{\begin{eqnarray*}}
\newcommand{\eeqno}{\end{eqnarray*}}
\newcommand{\beqla}[1] {\begin {eqnarray}\label{#1}}
\def\eeq {\end {eqnarray}}
\newcommand{\beq}{\begin {eqnarray}}

\newcommand{\real}{{\bf R}}

\newcommand{\integer}{{\bf Z}}

\newcommand{\complex}{{\bf C}}

\newcommand{\smooth}{C^\infty}

\newcommand{\F}{{\mathcal F}\, }
\newcommand{\FI}{{\mathcal F}^{-1}}

\newcommand{\sgn}{{\rm sgn}\,}

\newcommand{\torus}{{\bf T}}

\newcommand{\BMATRIX}[1]{\begin{bmatrix}#1\end{bmatrix}}

\newcommand{\ftn}{\mathcal{F}}

\newcommand{\E}{\mathbb{E}}
\newcommand{\vare}{\varepsilon}
\newcommand{\al}{\alpha}
\newcommand{\mass}{\mathbb{P}}
\DeclareMathOperator{\htt}{{\mathcal H}}

\DeclareMathOperator{\ba}{{\mathcal B\!\mathcal A}}

\DeclareMathOperator{\umd}{{\rm UMD}}
\newcommand{\prob}{{\mathbf P}}
\newcommand{\expec}{{{\mathbf E}\,}}

\DeclareMathOperator{\bsym}{{A_{\rm s}}}
\DeclareMathOperator{\basym}{{A_{\rm as}}}
\newcommand{\bsymm}[1]{{A_{{\rm s},#1}}}

\newcommand{\SYM}{{\scriptstyle{\BMATRIX{1&0\\ 0&-1}}}}
\newcommand{\ASYM}{{\scriptstyle\BMATRIX{0&-1\\ 1&0}}}
\newcommand{\ito}[3]{\| #1 |(#2 ,W)\|_{#3}}
\newcommand{\dito}[3]{\| #1 |(#2 ,g)\|_{#3}}

\newcommand{\sphe}{{\mathbf s\mathbf h}}

\newcommand{\ds}{\displaystyle}
        
        \def\halmos{{\ \vbox{\hrule\hbox{\vrule
height1.3ex\hskip0.8ex\vrule}\hrule}}
            \par\medskip}
        
        \newcommand{\refeq}[1]{(\ref{#1})}

\begin{document}
\title[On singular integral and martingale transforms]
{On singular integral and martingale transforms}

\author{Stefan Geiss}
\address{University of Jyv\"askyl\"a, Department of Mathematics and Statistics,
         P.O. Box 35 (MaD), FIN-40014 University of Jyv\"askyl\"a, Finland}
\email{geiss@maths.jyu.fi}

\author{Stephen Montgomery-Smith}
\address{University of Missouri, Mathematics Department,
         Columbia, MO 65211 USA}
\email{stephen@math.missouri.edu}

\author{\break{Eero Saksman}}
\address{University of Helsinki, Department of Mathematics and Statistics,
         P.O. Box 68, FIN-00014 University of Helsinki, Finland}
\email{eero.saksman@helsinki.fi} 


\thanks{The first and the last author are supported by the Project \#110599
        of the Academy of Finland.}

%%%%%%%%%%%%%%%%%%%%%%%%%%%%%%%%%%%%%% %%%%%%%%%%%%%%%%%%%%%%%%%%%%%%%%%%%%%%%%%%

\begin{abstract}  
Linear equivalences of norms of vector-valued singular integral operators and 
vector-valued martingale transforms are studied. In particular, 
it is shown that the UMD-constant of a Banach space $X$ equals 
the norm of the real (or the imaginary) part of the Beurling-Ahlfors singular 
integral operator, acting on $L^p_X(\real^2)$ with $p\in (1,\infty).$ 
Moreover, replacing equality by a linear equivalence, this is found to be the 
typical property of even multipliers.  A corresponding result for odd multipliers 
and the Hilbert transform is given. As a corollary we obtain that the norm of the 
real part of the Beurling-Ahlfors operator equals $p^*-1$ with 
$p^*:= \max \{p, (p/(p-1))\}$, where the novelty is the lower bound. 
\end{abstract} 

\maketitle

%%%%%%%%%%%%%%%%%%%%%%%%%%%%%%%%%%%%%%%%%%%%%%%%%%%%%%%%%%%%%%%%%%%%%%%%%%%%%%%%%

\section{Introduction}\label{se:intro}
A Banach space $X$ is said to be a UMD-space  
\footnote{UMD stands for {\em unconditional martingale differences}.}
provided that for all (equivalently, for some)
$p\in (1,\infty)$ there is a constant  $c_p>0$ such that
\[   \sup_{\alpha_k\in \{\pm 1\}} 
              \left \| \sum_{k=1}^n \alpha_k  D_k \right \|_{L_X^p}
   \le c_p \, \left \| \sum_{k=1}^n           D_k \right \|_{L_X^p} \]
for all $n\geq 1$ and all $X$-valued martingale difference sequences
$(D_k)_{k=1}^n$. As UMD-constant one usually takes $\umd_p(X):=\inf c_p$. We refer to
\cite{Burk3} and the references therein for an overview 
 about the UMD-property.
It is known that in the above definition the arbitrary martingale differences
can be replaced by Walsh-Paley martingale differences and one gets the same 
constant (see e.g. \cite[p. 12]{Burk4}
and \cite{Maurey}; the definition of Walsh-Paley martingales is recalled in Section \ref{se:definitions} below). The UMD-property was first investigated by Burkholder \cite{Burk1},
who gave a geometric characterization for Banach spaces with the 
UMD-property.  
Together with McConnell, Burkholder established in \cite{Burk2} that the 
Hilbert transform 
\beqla{eq:1.1}
\htt f(x):={\frac{1}{\pi}}\int_\real {\frac{f(y)}{x-y}}\, dy
\eeq
is bounded on $L^p_X(R)$, $p\in (1,\infty)$, provided that $X$ has the UMD-property. 
A converse result was proved soon after by Bourgain \cite{Bou},
who showed that the boundedness of the Hilbert transform on $L^p_X$ for
some $p\in (1,\infty )$ implies the UMD-property for $X.$ It is
also known from \cite{Bou2} that the UMD-property implies the boundedness of 
all invariant singular integrals or (more generally)  standard  multiplier 
operators under some regularity assumptions. The importance of the 
UMD-property, especially 
in connection with  PDE:s, is further evidenced by recent new results
on operator valued singular integrals (R-boundedness) \cite{Weis} and other
developments.

It is natural to ask for the quantitative equivalence of the UMD-property and the 
boundedness of the vector-valued Hilbert transform. The proofs in \cite{Burk2} 
and \cite{Bou} yield that there is a constant $C>0$ such that 
\beqla{eq:1.10}
     {\frac{1}{C}}(\umd_2 (X))^{1/2}\leq \| \htt\|_{L^2_X(\real )\to L^2_X(\real )} 
\leq C(\umd_2 (X))^2.
\eeq
The curious feature above is the  {\it quadratic equivalence} of the norms  in 
contrast to the linear dependence one would expect. A better than a quadratic 
equivalence obtained from alternative proofs is not known to the authors.

The previous discussion raises the question whether there is a linear equivalence 
in \refeq{eq:1.10}. We do not know whether this is true or not.
%but we  suspect that the answer is negative. 
However, in this paper we establish that the answer is positive
if $\htt$ is replaced by the Beurling-Ahlfors transform $\ba$:
\medskip

\thm{th:main}
For $p\in (1,\infty )$ and a real Banach space $X$ one has that
\beqla{eq:1.20}
     \hspace*{1.5em}
     \umd_p (X) = \| {\rm Re}(\ba)\|_{L^p_X(\real^2 )\to L^p_X(\real^2 )} 
                = \| {\rm Im}(\ba)\|_{L^p_X(\real^2 )\to L^p_X(\real^2 )} 
\eeq
with
\beqla{eq:1.21}
\ba f(z):=-{\frac{1}{\pi}}\int_{\complex} {\frac{f(w)}{(z-w)^2}}\, dm_2(w),
\eeq
where $m_2$ is the two-dimensional Lebesgue measure on the complex plane $\complex$, that  is identified 
with $\real^2.$
\ethm
\medskip

Note that  ${\rm Re}(\ba) = - Id - 2 R_1^2$ and ${\rm Im}(\ba)=2 R_1R_2$
where $R_1$ and $R_2$ are the first and second Riesz transform, respectively.
Equality (\ref{eq:1.20}) carries some  new information for norm estimates even in 
the scalar case, see Corollary \ref{cor:3.1} below and the remark after it.  

The operator $\ba$ is sometimes called the two-dimensional Hilbert transform. 
It plays a fundamental role in the theory of quasi-conformal maps and in the 
theory of elliptic equations in the plane. Quite recently,
in connection with the well-known Iwaniec conjecture,
there has been many works devoted to  the probabilistic approach to
estimate the (scalar) $L^p$-norm of the Beurling-Ahlfors operator, see e.g.
\cite{Ban1,Volberg,Ban2}.

If we replace the equality in (\ref{eq:1.20}) by a linear equivalence
with multiplicative constants, then the obtained property is shared by an extensive class 
of operators
corresponding to homogeneous multipliers. In this context the crucial difference
between the Hilbert transform and the Beurling-Ahlfors operator is the parity of 
their integral kernels: our main Theorem \ref{th:3.1} shows 
that {\it one may replace the real and imaginary part of
$\ba$ in {\rm \refeq{eq:1.20}} by
 any smooth, homogeneous of order zero, and even
Fourier  multiplier operator} if we allow multiplicative constants.
On the other side, Theorem \ref{th:4.1} in turn relates  odd multipliers
to the Hilbert transform.

The proof of Theorem \ref{th:main} is based on a modification of an argument 
of Bourgain and the already standard representation of suitable combinations 
of Riesz singular integrals as certain transforms of Ito-integrals,
see \cite{GunVar,Ban1}. In order to define a (wider) class of such 
transforms, let $W=((W_t^1,...,W_t^d))_{t\ge 0}$ be a standard $d$-dimensional
Brownian motion. For a Banach space $X$, $p\in (1,\infty)$, and a real  
$d\times d$ matrix $A$ we let $\ito{Id_X}{A}{p} := \inf c$, where the infimum 
is taken over all $c>0$ such that
\[     \Big\| \int_0^T U\cdot  d(A W_t)
       \Big \|_{L^p_X}
   \le c  \Big \| \int_0^T U\cdot  d(W_t)
          \Big \|_{L^p_X} \]
for all $T\ge 0$ and  certain $d$-tuples $U$ of $X$-valued processes  
taking values in a finite-dimensional subspace of $X$ 
(see Section \ref{se:definitions}).
As a byproduct of our proofs we obtain in Theorem \ref{th:3.1}
that, if the matrix $A$ is
symmetric and not a constant multiple of the identity matrix,
then there is a constant $C=C(A)>0$, independent of the Banach space $X$,
such that
\beqla{eq:1.30}
     {\frac{1}{C}}\umd_p (X)
\leq \ito{Id_X}{A}{p}
\leq C \umd_p (X).
\eeq

In the proof of Theorem \ref{th:main} (in fact, Theorem \ref{th:3.1})
we do not employ  harmonic extensions of functions 
to the upper half plane as it was done
originally by Gundy and Varopoulos \cite{GunVar} in their stochastic proof of
the $L^p$-boundedness 
of the scalar Riesz transforms. Instead, we follow Ba{\~n}uelos and Mendez-Hernandez
\cite{Ban2} and use the space-time Brownian motion, which corresponds to the heat 
extension of functions. Earlier the use of heat extensions in norms estimates for singular integrals 
was initiated by Nazarov, Petermichl and Volberg in \cite{Petermichel-Volberg} and \cite{Volberg}.
The use of the space-time Brownian motion in combination with our modified version of  Bourgain 
techniques is one reason that enables us to obtain the equalities (\ref{eq:1.20}). In particular, for all $d\geq 2$ it holds that
\begin{equation}\label{eq:1.32}
\umd_p(X) = \left \|  Id+2R_1^2 \right \|_{L^p_X(\real^d )\to L^p_X(\real^d )}
          = \ito{Id_X}{A_s}{p}
\end{equation}
with
\[ A_s := \SYM. \]

As one specializes to the case $X=\real$ in Theorem \ref{th:main} one obtains
\medskip

\begin{corollary}\label{cor:nice}
For $p\in (1,\infty)$ one has that
\[   \|{\rm Re}(\ba) \|_{L^p(\real^2)\to L^p(\real^2)}
   = p^*-1
   \qquad \mbox{where}\;\; 
   p^*:=\max \left (p,\frac{p}{p-1}\right ). \]
\end{corollary}
\noindent This follows  because Burkholder 
(see \cite{Burk5},
     \cite[Theorem 14]{Burk3},
     \cite[p.12]{Burk4}) 
established that 
$\umd_p(\real )=p^*-1$. One should also observe here, that it is well-known that the $L^p-$norm of
scalar multipliers that preserve real-valued functions does not depend on whether one look at
real or complex $L^p$-spaces. The upper bound in Corollary \ref{cor:nice} is due to
Volberg and Nazarov \cite[Theorem 3.1]{Volberg} who applied Bellman function and heat kernel techniques
in the proof (a space-time Brownian motion approach 
to  \cite[Theorem 3.1]{Volberg} was subsequently 
given in \cite{Ban2}).
The real novelty in our Corollary \ref{cor:nice} is the lower bound. As far as we
know, this is the first instance where the exact value of the  norm of a 'second order' Riesz transform
(a non-trivial expression involving products of Riesz transforms) has been determined.
Moreover,  Corollary \ref{cor:nice} has obvious interest in connection with the
 famous Iwaniec conjecture, which claims that the norm of the {\it full} complex Beurling-Ahlfors transform equals $p^* -1.$

The paper is organized as follows:
Section \ref{se:definitions} recalls necessary notation and contains
preparatory material. Our main result, Theorem \ref{th:3.1}, which
contains Theorem \ref{th:main}, \refeq{eq:1.30}, and \refeq{eq:1.32} as  
special cases, is formulated and proved in Section \ref{se:results}.  
Section \ref{se:hilbert} treats the case of odd kernels and
anti-symmetric matrix transforms of stochastic integrals. 

The results in Sections \ref{se:results} and \ref{se:hilbert} are formulated for
tensor product operators $T_m\otimes S$, where $S: X\to Y$ is an arbitrary
operator between two Banach spaces $X$ and $Y$, instead for the setting of
identities $Id_X$ we have used so far in the introduction.
The motivation for this are certain connections to the geometry of Banach spaces
that are explained in Section \ref{se:remarks}.

We would like to thank Tuomas Hyt\"onen for the careful reading of the
manuscript. In his recent paper \cite{Hyt} he obtained results about the
linear equivalence of norms of vector valued spectral multipliers 
and the UMD-constant, which are in the spirit of the results in this 
paper (but, as pointed out in \cite{Hyt}, the basic singular integral operators
like $\ba$ or $\htt$ do not fall into the setting and scope of \cite{Hyt}).

%%%%%%%%%%%%%%%%%%%%%%%%%%%%%%%%%%%%%%%%%%%%%%%%%%%%%%%%%%%%%%%%%%%%%%%%%%%%%%%%%

\section{Definitions and preliminary results}
\label{se:definitions}

We shall denote by $\{ e_1,\ldots ,e_d\}$ the unit vectors of $\real^d$,
by $|x|$ the euclidean norm of $x\in\real^d$, and will use 
$B(x,\delta):= \left \{ y\in\real^d : |x-y|\le\delta \right \}$ 
for $\delta\ge 0$. Moreover, the set of real $m\times n$ matrices is 
denoted by $M(m,n).$
\medskip 

\paragraph{{\bf Vector-valued operators}.}  In this paper $L(X,Y)$ stands 
for the linear and bounded operators between two real Banach spaces $X$ and $Y$, 
where $L(X):=L(X,X)$. Given a $\sigma$-finite measure space
$(M,\mu)$ and $p\in [1,\infty)$, the space of Bochner integrable random variables 
$L^p_X(M)=L^p_X(M,\mu )$ consists of all 
strongly measurable functions  $f:M \rightarrow X$ such that there is a
separable subspace $X_0 \subseteq X$ with $f(M) \subseteq X_0$ and 
$   \| f\|_{L^p_X} 
 := \left ( \int_M \| f \|^p d\mu \right )^{{1}/{p}}
  < \infty$.
Given operators  $S\in L(X,Y)$ and $T\in L(L^p(M,\mu))$
the tensor product $T\otimes S$  can be defined 
through its action on simple functions, i.e.
\[     (T\otimes S)\left ( \sum_{k=1}^n x_k\chi_{E_k}\right) 
    := \sum_{k=1}^n Sx_kT(\chi_{E_k}). \]
In case that the quantity
\begin{multline*}
      \|T \otimes S:L^p_X(M,\mu) \rightarrow L^p_Y(M,\mu)\| \\
   := \sup \left \{ \|  (T\otimes S) (F) \|_{L^p_Y} : \| F\|_{L^p_X} \le 1,
                       F \;\mbox{ simple function}\right \}
\end{multline*}
is finite, the operator $T\otimes S$ extends to a bounded linear operator from
$L_X^p(M,\mu)$ into $L_Y^p(M,\mu)$. As particular operators $T$ we use 
multipliers. 
The usage of operators $S$ instead of identities ${\rm Id}:X\to X$ of Banach 
spaces $X$ might be seen as somehow artificial at this point. 
As already mentioned, the motivation for this slightly more general setting can be found in 
Section \ref{se:remarks}, but the reader may, if she or he so wishes, replace the 
operator $S$ by an identity in what follows. 
\bigskip

\paragraph{{\bf Multipliers}.}
A bounded complex valued function $m\in \smooth (\real^d\setminus \{ 0\})$,
$d\ge 1$, is called (smooth) multiplier. A multiplier $m$ is homogeneous 
(of order zero) if $m(\lambda \xi )=m(\xi )$ for $\xi\in\real^d\setminus \{ 0\}$ 
and $\lambda >0$.
In this paper the term multiplier always refers to smooth and homogeneous multipliers. 
The multiplier $m$ is called even provided that $m(\xi )=m(-\xi )$ for 
$\xi\not =0$ and odd if $m(\xi )=-m(-\xi )$.
The operator $T_m:L^2(\real^d) \rightarrow L^2(\real^d)$ associated to $m$ 
is given by
\[ T_mf := \FI (m\F f), \]
where $\F$ stands for the Fourier transform
\[ (\F f)(\xi):=\int_{\real^d} \exp(-i\langle \xi,x \rangle) f(x) dx. \]
It is easy to check that an even and real multiplier maps real-valued 
functions to real-valued functions, and the same is true for odd and purely
imaginary multipliers. Consequently, to formulate our results for real Banach 
spaces, from now on we use the standing 
\begin{assumption}\label{assumption:m}\rm
All multipliers are even and real, or odd and purely imaginary.
\end{assumption}

Let $A\in M(d,d)$ be invertible. By applying the simple identity
\[ T_{m\circ A}f(x) = (T_m(f\circ A^T))\left ((A^{T})^{-1}x\right ) \]
we deduce that composing a multiplier with a linear invertible 
map one does not change its norm, i.e.

\begin{equation}
\label{eq:2.10}
  \|{\rm T_{\it m\circ A}} \otimes S:L^p_X(\real^d) \rightarrow L^p_Y(\real^d)\|
= \|{\rm T_{\it m}       } \otimes S:L^p_X(\real^d) \rightarrow L^p_Y(\real^d)\|.
\end{equation}

For a multiplier $m$ on $\real^d$ as above there is a corresponding 
discrete multiplier $\widetilde m$ that acts on functions defined on the 
$d$-dimensional torus $\torus^d:=(-\pi,\pi ]^d$: for  
a finite trigonometric polynomial $f$ we let
\[     (T_{\widetilde{m}}f)(\theta)
   := \sum_{k\in\integer^d}\hat  f(k)e^{i \langle k,\theta\rangle}m(k),
\]
where 
$\hat f(k) := (1/2\pi)^d \int_{\torus^d}e^{-i \langle k,\theta \rangle}
             f(\theta)d\theta $ 
and $m(0):=\omega_{d-1}^{-1}\int_{S^{d-1}}m(x)dx$
is the average over the boundary of the euclidean unit ball
(remember that $m$ is homogeneous of order zero). It follows by 
Assumption \ref{assumption:m} that in the above definition 
$T_{\widetilde{m}}f$ is real whenever $f$ is real. 
In the sequel it will be important that the norms of the corresponding 
multipliers are equal which will be stated in 
Lemma \ref{le:multiplier-torus-rn} below.
\smallskip
\lem{le:multiplier-torus-rn}
Let $m$ be a  smooth and homogeneous multiplier on $\real^d$, $d\ge 1$.
Then, for any $p\in (1,\infty)$ and $S\in L(X,Y)$, it holds that
\begin{eqnarray*}
      \| T_m\otimes S : L^p_X(\real^d )\to L^p_Y(\real^d ) \| 
& = & \| T_{\widetilde m}\otimes S : L^p_X(\torus^d )\to L^p_Y(\torus^d )\|\\
& = & \| T_{\widetilde m}\otimes S : 
            L^p_{X,0}(\torus^d )\to L^p_Y(\torus^d )\| , 
\end{eqnarray*}
where $L^p_{X,0}(\torus^d)$ stands for the functions in $L^p_{X}(\torus^d)$ 
of mean zero.
\elem
\smallskip

\begin{proof}
The first equality in the scalar case is essentially due to K. de Leeuw 
\cite{deLeeuw}, see also \cite{CoifmanWeiss}.
The proof in the monograph \cite[pp. 221--223]{Grafakos} can
be easily verified to carry over to the case considered here to yield
the estimate
\[
\| T_{\widetilde m}\otimes S : L^p_X(\torus^d )\to L^p_Y(\torus^d )\|
\le \| T_m\otimes S : L^p_X(\real^d )\to L^p_Y(\real^d ) \|.
\]
The lemma follows as soon as one has 
\begin{equation*}%\label{eq:2.1001}
\| T_m\otimes S : L^p_X(\real^d )\to L^p_Y(\real^d ) \| \leq
\| T_{\widetilde m}\otimes S : L^p_{X,0}(\torus^d )\to L^p_Y(\torus^d )\|.
\end{equation*}
Observe that when computing the norm of $T_m\otimes S$ we may restrict 
ourselves to functions of type $f=\sum_{k=1}^n f_k x_k$ where $x_k\in X$ 
and the $f_k\in C_0^\infty(\real^d)$ have integral zero
because these functions $f$ are dense in $L^p_{X}(\real^d)$. 
Then we can also follow the second part of the proof of the transference
principle in \cite{Grafakos} on pp. 223--225. One also verifies, that in the above proofs
one may restrict to real valued scalar functions $f_k$, and as  our multipliers
preserve real-valued functions, the result remains true also for real Banach spaces $X,Y.$
\halmos
\end{proof}
\bigskip

\paragraph{{\bf Hilbert transform}.} 
The Hilbert transform was defined in \refeq{eq:1.1}. It corresponds to
the multiplier $m(\xi)= -i\,\sgn (\xi )$, and maps, by definition,
real valued functions to real valued functions. The corresponding
discrete multiplier operator $T_{\widetilde m}$ is the well known
conjugation operator 
$\widetilde\htt :L_p(\torus)\rightarrow L_p(\torus)$, $p\in (1,\infty)$,
which can be also defined through its action on the trigonometric
polynomials
\[ (\widetilde\htt \sin (k\cdot))(\theta) :=   - \cos (k \theta)
   \sptext{1}{and}{1}
   (\htt \cos (k\cdot))(\theta) :=   \sin (k\theta) \] 
for $k=1,2,...$ and $\widetilde\htt 1:\equiv 0$.
\medskip

\paragraph{{\bf Beurling-Ahlfors transform and Riesz transforms}.} 
The $k$-th Riesz-trans\-form $R_k$, $k=1,\ldots ,d$ is the multiplier operator on $\real^d$
corresponding to the multiplier $\xi_k/(i|\xi | )$. 
The Beurling-Ahlfors operator, defined through
\refeq{eq:1.21} in the introduction, corresponds to the multiplier
$m_{BA}(\xi ):=\ds  {\frac{\xi_1-i\xi_2}{\xi_1+i\xi_2}}$. It follows that
\[ \ba = R_2^2-R_1^2+2iR_1R_2
       =: {\rm Re}(\ba) + i {\rm Im}(\ba). \]
For our real valued setting  we consider ${\rm Re}(\ba)$ and ${\rm Im}(\ba)$
separately and check by \refeq{eq:2.10} that a rotation of the 
coordinates by the angle $\pi/4$ transforms ${\rm Re}(\ba)$ into 
${\rm Im}(\ba)$ and their norms coincide.
\smallskip

\paragraph{{\bf Martingale transforms}.} 
Given independent Bernoulli random variables $\vare_1,\vare_2,\ldots$, i.e.
$\prob(\vare_k=\pm 1)=1/2$, and maps $d_k: \real^{k-1}\to X$, where $X$ is a Banach space
and $d_1$ is constant, a sequence $(\vare_k d_k(\vare_1,...,\vare_{k-1}))_{k\in I}$
with $I=\{1,...,n\}$ or $I=\{1,2,...\}$ is called Walsh-Paley martingale difference sequence.
Given $S\in L(X,Y)$ and $p\in (1,\infty)$, we define
$\umd_p(S) := \inf c$, where the infimum is taken over all $c>0$ such that
\beqla{eq:2.79}       
         \left \| \sum_{k=1}^n \alpha_k S \vare_k d_k \right \|_{L^p_Y}
   \le c \left \| \sum_{k=1}^n \vare_k d_k \right \|_{L^p_X} 
\eeq
for all $X$-valued Walsh-Paley martingale difference sequences $(\vare_k d_k)_{k=1}^n$, 
all $\alpha_1,...,\alpha_n\in\real$ with $|\alpha_k|\le 1$, and  
$n\ge 1$. It is well-known that by an easy extreme point argument  
the condition $|\alpha_k| \le 1$ can be replaced by $\alpha_k\in\{-1,1 \}$ and one gets
the same constant $\umd_p(S)$. 
The reader is also referred to \cite{Burk3} and the references therein for a more
general overview about UMD-spaces.
\medskip
 
\paragraph{{\bf Transforms for stochastic integrals}.}
We recall the definition given in the introduction.
Let $W=((W_t^1,...,W_t^d))_{t\ge 0}$ be a $d$-dimensional
standard Brownian motion with continuous paths for all $\omega\in\Omega$
and $W_0\equiv 0$, defined on a probability space $(\Omega,\ftn,\mass)$,
where $\ftn$ is the completion of $\sigma(W_t: t\ge 0)$ and 
$(\ftn_t)_{t\ge 0}$ the augmentation of the natural filtration
of $W$. Let $S\in L(X,Y)$, $p\in (1,\infty)$, and $A=[a_{kl}]\in M(d,d)$. 
Then $\ito{S}{A}{p} := \inf c$, such that
\begin{equation}\label{eqn:definition_rotation_integrand}
       \Big \| \sum_{k=1}^d \sum_{l=1}^d a_{kl} \int_0^T SU_t^k d W_t^l 
       \Big \|_{L^p_Y}
   \le c  \Big \| \sum_{k=1}^d \int_0^T U_t^k d W_t^k
          \Big \|_{L^p_X}
\end{equation}
for all $T\ge 0$ and $(\mathcal{F}_t)_{t\ge 0}$-adapted left-continuous processes 
of Radon random variables $U_t^k:\Omega\rightarrow X$ which have right-hand side limits,
take values in a finite-dimensional subspace of $X$, satisfy
\[ \int_0^T \E\| U_t^k\|_X^2 dt < \infty \]
for all $T\ge 0$ and $k=1,...,d$, and such that the right hand side of 
(\ref{eqn:definition_rotation_integrand}) is finite.
To shorten the notation we also use in the sequel 
$\int_0^T (SU_t^k)_{k=1}^d \cdot d(AW)_t$ and
$\int_0^T (U_t^k)_{k=1}^d \cdot dW_t$, respectively, 
for the expressions inside the norms.

\bigskip
\paragraph{{\bf Some minor notation}.} Given $A,B\ge 0$ and $c>0$, the notation
$A\sim_c B$ stands for $A/c\le B \le cA$. If the dependence of $c$ on the
extra quantities involved is clear, we sometimes simply write $A\sim B.$

%%%%%%%%%%%%%%%%%%%%%%%%%%%%%%%%%%%%%%%%%%%%%%%%%%%%%%%%%%%%%%%%%%%%%%%%%%%%%%%

\section{The main result}\label{se:results}

It will be convenient to have a special notation for particular matrices and multipliers. 
Thus, we denote 
\begin{equation}\label{eqn:matrices}
 \hspace*{.4em}    
     \bsym := \SYM, 
     \quad 
 \bsymm{d} := {\scriptstyle\BMATRIX{1&0&0& \hdots&0\\ 0&-1&0&\hdots&0\\
              0&0&-1& \hdots&0\\\vdots &\vdots&\vdots&\ddots&\vdots\\0&0&0&0&-1}}, 
     \quad 
    \basym := \ASYM .
\end{equation}
We shall also define the special multiplier $m_0$ on 
$\real^d$, where $d\geq 2$, by
\[
m_0 (\xi ) := 2{\frac{\xi_1^2}{|\xi |^2}} - 1, \quad \mbox{equivalently}\quad 
T_{m_0}    := - ({\rm Id} + 2R_1^2).
\]

We recall once more that, for simplicity, all Banach spaces are assumed to be real.
Our  main result is the following:
\bigskip

\thm{th:3.1}
Assume that $m\in\smooth (\real^d\setminus \{0\})$, $d\ge 2$,  is a real, homogeneous, and 
even multiplier that is not identically constant, and $A\in M(d,d)$ is a real symmetric
matrix that is not a multiple 
of the identity matrix. Let $p\in (1,\infty ).$ Then there is
a constant $C=C(m,A)$ such that for every pair of Banach spaces $X$ and $Y$ and for every 
operator $S\in L(X,Y)$ one has that
\begin{multline*}
      \|  S| ({\bsym} ,W)\|_p
    =  {\rm UMD}_p(S) 
    = \| T_{m_0} \otimes S : L^p_X(\real^d)\to L^p_Y(\real^d) \| \\
   \sim_C \| T_m\otimes S : L^p_X(\real^d)\to L^p_Y(\real^d) \|
   \sim_C  \|  S| (A,W)\|_p.
\end{multline*} 
\ethm 
\bigskip
Exploiting $\umd_p(\real )=p-1$ for $p\in [2,\infty)$ (see \cite{Burk5},
\cite[Theorem 14]{Burk3}) and 
using that, for $d=2$, the multiplier $m_0$  corresponds to $R_2^2-R_1^2$ 
and $2R_1R_2$ can be obtained by a rotation of $m_0$, 
Theorem \ref{th:3.1} implies 
\bigskip

\begin{corollary}\label{cor:3.1} 
For $d=2$, $p\in (1,\infty)$, and $S\in L(X,Y)$ one has
\begin{multline*}
     \| {\rm Re}(\ba) \otimes S : L^p_X(\real^2)\to L^p_Y(\real^2)\| \\
   = \| {\rm Im}(\ba) \otimes S: L^p_X(\real^2)\to L^p_Y(\real^2)\| 
   = {\rm UMD}_p(S).
\end{multline*}
In particular, for $p\in [2,\infty)$,
\[   \| {\rm Re}(\ba): L^p(\real^2)\to L^p(\real^2)\| 
   = \| {\rm Im}(\ba): L^p(\real^2)\to L^p(\real^2)\| 
   = p-1. \]
\end{corollary}
\bigskip 

% In the second part of Corollary \ref{cor:3.1} 
% the upper bound $p-1$ was known, see \cite{Ban2}.
% The lower bound is of interest in connection to the Iwaniec conjecture 
% (see e.g. \cite{Volberg,Ban2,Baernstein-MS})
% as it follows  that  for $p \in [2,\infty)$ the norm of just the real 
% part of $\ba$ acting on real valued functions gives already
% $p-1$.
% \bigskip

We break the proof of Theorem \ref{th:3.1} into a series of auxiliary results,
some of which are basically known, and some of which may have independent interest.
The actual proof of Theorem \ref{th:3.1} is given at the end of this
section. Our first step is to show that non-trivial even multipliers dominate
the UMD-constant linearly. To this end we need a generalization 
of a Lemma due to
Bourgain  \cite[Lemma 1]{Bou}. For the rest of this section we assume that
$X$ and $Y$ are Banach spaces and  $S\in L(X,Y)$ is a fixed operator.
\bigskip

\lem{le:Bourgain_new}
Let  $p\in (1,\infty)$, $Q:=\torus^d$, and assume that the 
multiplier $m$ satisfies Assumption \ref{assumption:m}. For $k\ge 1$
let $E_k$ be the closure in $L^p_X(Q^k)$ of the finite real trigonometric 
polynomials
\[   \Phi_k(\theta_1,...,\theta_k) 
   = \sum_{p=1}^P  \Phi_k^p(\theta_1,...,\theta_k) x_p \]
with $x_k\in X$ and
\[   \Phi_k^p
   = \sum_{\ell_1\in\integer^d} \cdots \sum_{\ell_k\in\integer^d}
     e^{i \langle \ell_1,\theta_1 \rangle} \cdots
     e^{i \langle \ell_k,\theta_k \rangle} \alpha^p_{\ell_1,...,\ell_k} 
   \sptext{1}{where}{1}
   {\rm Im}(\Phi_k^p)\equiv 0 \]
where only finitely many of the $\alpha^p_{\ell_1,...,\ell_k}\in \complex$
are non-zero and 
$\alpha^p_{\ell_1,...,\ell_k} = 0$ whenever $\ell_k=0$ (so that
$\int_Q \Phi_k^p(\theta_1,...,\theta_k) d \theta_k = 0$).
Let
$T_{\widetilde{m}}^k: L^p(Q^k) \to L^p(Q^k)$ be given by
\begin{multline*}
     (T_{\widetilde{m}}^k \Phi_k)(\theta_1,...,\theta_k) \\
   = \sum_{p=1}^P \left (
     \sum_{\ell_1\in\integer^d} \cdots \sum_{\ell_k\in\integer^d}
     m(\ell_k)
     e^{i \langle \ell_1,\theta_1 \rangle} \cdots
     e^{i \langle \ell_k,\theta_k \rangle} \alpha^p_{\ell_1,...,\ell_k}
                   \right ) x_p.
\end{multline*}
for $\ell_1,...,\ell_k\in\integer^d$. 
Then one has that
\begin{multline*}%\label{eq:3.20_new}
     \left \| \sum_{k=1}^n 
     ((T_{\widetilde m}^k \otimes S) \Phi_k)(\theta_1,\ldots,\theta_k)
     \right \|_{L^p_Y(Q^n)} \\
\leq \| T_{\widetilde m}\otimes S: {L^p_X(\torus^d)\to L^p_Y(\torus^d)}  \|
     \left  \| \sum_{k=1}^n\Phi_k (\theta_1,\ldots,\theta_k)\right \|_{L^p_X(Q^n)}
\end{multline*}
for $\Phi_1\in E_1$,..., $\Phi_n\in E_n$.
\elem
\smallskip

\begin{proof} 
It is sufficient to prove the inequality for finite real trigonometric 
polynomials 
$\Phi_1(\theta_1),...,\Phi_n(\theta_1,...,\theta_n)$.
Let $A\geq 1$ be an integer, $\eta\in Q$ be an auxiliary variable, denote 
by $T_{{\widetilde m},\eta}$ the application of $T_{\widetilde m}$ with respect to the 
variable $\eta$, and consider the difference
\begin{eqnarray*}
      D_k^A(\theta_1,...,\theta_k,\eta)
&:= & ((T_{{\widetilde m},\eta} \otimes S) \Phi_k(\theta_1+A\cdot,\ldots ,\theta_k+A^k\cdot))
      (\eta) \\
&   & -  ((T_{\widetilde m}^k \otimes S) \Phi_k) 
      (\theta_1+A\eta,\ldots ,\theta_k+A^k\eta ).
\end{eqnarray*}
Note, that in the first term on the right-hand side we apply the multiplier 
$T_{\widetilde m}$ to a function, where
$\theta_1,...,\theta_k,A$ act as parameters, whereas in the second term 
we apply the multiplier $T_{\widetilde m}^k$ to $\Phi_k$ itself.
If we show that for our fixed trigonometric polynomial $\Phi_k$ there is the bound
\begin{equation}\label{eq:3.40}
\sup_{A\ge 1} A \| D_k^A\|_{L_Y^\infty}  < \infty,
\end{equation}
then the proof is completed by, firstly, observing  that  
\begin{multline*}
    \left \| \sum_{k=1}^n 
    ((T_{{\widetilde m},\eta} \otimes S) \Phi_k(\theta_1+A\cdot,\ldots ,\theta_k+A^k\cdot))
    (\eta) 
    \right\|_{L_Y^p(Q,d\eta)} \\
\le \| T_{\widetilde m}\otimes S: {L^p_X(Q)\to L^p_Y(Q)}\| 
     \left \| \sum_{k=1}^n 
     \Phi_k (\theta_1+A\eta,\ldots ,\theta_k+A^k\eta )
     \right\|_{L_X^p(Q,d\eta)}.
\end{multline*}
Secondly, one replaces  the left-hand side by
\[ \left \| \sum_{k=1}^n 
   ((T_{\widetilde m}^k \otimes S) \Phi_k) (\theta_1+A\eta,\ldots ,\theta_k+A^k \eta )
   \right\|_{L_Y^p(Q,d\eta)} \]
with the corresponding correction terms in $A$,
integrates with respect to the $\theta$'s, 
applies Fubini's theorem so that the variable  $\eta$ is removed, and 
sends $A$ to infinity.

In order to prove \refeq{eq:3.40} we observe that 
$A \| D_k^A(\theta_1,...,\theta_k,\eta) \|_Y$ is upper bounded
by a finite number of terms of form
\begin{multline*}
      A \bigg | 
      T_{{\widetilde m},\eta} \left (
      e^{i\langle \ell_1,\theta_1+A   \eta\rangle} \cdots  
      e^{i\langle \ell_k,\theta_k+A^k \eta\rangle} \right ) \\
      - m(\ell_k) 
      e^{i\langle \ell_1,\theta_1+A   \eta\rangle} \cdots  
      e^{i\langle \ell_k,\theta_k+A^k \eta\rangle}
      \bigg | |\alpha^p_{\ell_1,...,\ell_k}| \|Sx_p\|_Y
\end{multline*}
with $\ell_k\not = 0$.
Finally we observe that 
\begin{eqnarray*}
&   & A \bigg | 
      T_{{\widetilde m},\eta} \left (
      e^{i\langle \ell_1,\theta_1+A   \eta\rangle} \cdots  
      e^{i\langle \ell_k,\theta_k+A^k \eta\rangle} \right ) 
      \\ 
&   & \hspace*{3em} - m(\ell_k) 
      e^{i\langle \ell_1,\theta_1+A   \eta\rangle} \cdots  
      e^{i\langle \ell_k,\theta_k+A^k \eta\rangle}
      \bigg | \\
& = & A |m(\ell_1 A         + \ldots + \ell_k A^k )- m(\ell_k)| 
      \\ \\
& = & A |m(\ell_1 A^{-k+1}  + \ldots  +\ell_{k-1} A^{-1} + \ell_k)- m(\ell_k)|
\end{eqnarray*}
which is bounded in $A$ since $m$ is differentiable at $\ell_k\not = 0$.
This yields \refeq{eq:3.40}.\halmos
\end{proof} 
\bigskip 
 
\prop{pr:UMD-multiplier}
Assume that $m\in C^\infty(\real^d\setminus \{ 0 \})$, $d\geq 2$, is a 
smooth, homogeneous, even, and non-constant multiplier. 
Let $\delta^+:=\max_{|\xi|=1}m(\xi )$ and $\delta^-:=\min_{|\xi|=1}m(\xi )$ 
so that $\delta^+-\delta^- >0.$  Then
\[ {\umd}_p(S)
   \leq {\frac{2}{\delta^+ -\delta^-}}
        \left( 1+ \frac{|\delta^+ + \delta^-|}{|\delta^+|+|\delta^-|}\right) 
        \| T_m\otimes S  
        : L^p_X (\real^d ) \to L^p_Y  (\real^d )\| . \]
In particular, in the case  $\max_{|\xi|=1}m(\xi )=-\min_{|\xi|=1}m(\xi )=1$
we have that
\[      {\umd}_p(S) 
   \leq \| T_m\otimes S : L^p_X(\real^d ) \to L^p_Y (\real^d) \|. \]
\eprop
\bigskip
 
\begin{proof}
By continuity and compactness there are 
$\xi^-,\xi^+\in\real^d$ of length one such that $m(\xi^-)=\delta^-$ 
and $m(\xi^+)=\delta^+$. Without loss of generality we may assume that 
$\xi^-=e_1$ and $\xi^+=e_2$, where $e_1$ and $e_2$ are the first two unit 
vectors in $\real^d$ (otherwise \refeq{eq:2.10} enables us to replace $m(\xi)$ by 
$m(A\xi )$ with suitably chosen $A$).
Define the functions $a^-,a^+\in L^\infty (\torus^d)$ by
$a^-(\theta):=\sgn (\theta_1)$ and  $a^+(\theta):=\sgn (\theta_2)$ for 
$\theta\in\torus^d$ so that
$T_{\widetilde{m}}a^- = \delta^- a^-$ and $T_{\widetilde{m}}a^+ =\delta^+ a^+$.
For independent Bernoulli random variables 
$\varepsilon_1,\varepsilon_2,\ldots$ we consider the $X$-valued 
Walsh-Paley martingale difference sequence 
\[ (\vare_k d_k(\vare_1,\ldots ,\vare_{k-1}))_{k=1}^n \]
and a sequence $(\al_k)_{k=1}^n$ with $\al_k \in \{ \delta^-,\delta^+ \}$. 
Define $\psi_k:=a^-$ if $\al_k=\delta^-$ and $\psi_k:=a^+$ if $\al_k=\delta^+$, and let
\[       \phi_k(\theta_1,\ldots ,\theta_{k-1})
   := d_k(\psi_1(\theta_1),\ldots,\psi_{k-1}(\theta_{k-1})). \]
Since $(\psi_1(\theta_1),\ldots ,\psi_n(\theta_n))$
and $(\vare_1,\ldots,\vare_n)$ have the same distribution 
(if we normalize the measure on $Q^n$) and 
$T_{\widetilde m}\psi_k = \alpha_k \psi_k$, Lemma \ref{le:Bourgain_new}
implies that 
\begin{multline*}
   \left \| \sum_{k=1}^n \alpha_k \vare_k S_k d_k(\vare_1,...,\vare_{k-1})
   \right \|_{L_Y^p} \\
   \le \| T_{\widetilde m}\otimes S: L^p_X(\torus^d)\to L^p_Y(\torus^d) \|
       \left \| \sum_{k=1}^n \vare_k d_k(\vare_1,...,\vare_{k-1})
       \right \|_{L_X^p}.
\end{multline*}
Let 
$A:=2/(\delta^+-\delta^-)$ and 
$B:=(\delta^+ +\delta^-)/(\delta^+-\delta^-)$ 
so that the new sequence $\beta_k := A \alpha_k - B$ satisfies 
$\beta_k=-1$ if $\alpha_k=\delta^-$ and  $\beta_k= 1$ if $\alpha_k=\delta^+$.
Then 
\begin{eqnarray*}
&   & \hspace*{-1em}
      \left \| \sum_{k=1}^n \beta_k \vare_k S_k d_k(\vare_1,...,\vare_{k-1})
      \right \|_{L_Y^p} \\
&\le& \!
      \left [ A  \| T_{\widetilde m}\otimes S: L^p_X(\torus^d)\to L^p_Y(\torus^d) \|
              + |B|  \| S \| \right ]
       \left \| \sum_{k=1}^n \vare_k d_k(\vare_1,...,\vare_{k-1})
       \right \|_{L_X^p}\!\!\!\!  .
\end{eqnarray*}
Because 
$      \| m \|_\infty \| S \| 
    =  \sup_{\ell\in\integer^d} |m(\ell)| \| S \| 
   \le \| T_{\widetilde m}\otimes S: L^p_X(\torus^d)\to L^p_Y(\torus^d) \|$
we end up with
\begin{eqnarray*}
      {\umd}_p(S)
&\le& \left [A + \frac{|B|}{\| m \|_\infty} \right ]
      \| T_{\widetilde m}\otimes S: L^p_X(\torus^d)\to L^p_Y(\torus^d) \| \\
& = & \left [A + \frac{|B|}{\| m \|_\infty} \right ]
      \| T_m\otimes S: {L^p_X(\real^d)\to L^p_Y(\real^d)} \| \\
&\le& \frac{2}{\delta^+ -\delta^-}
      \left( 1 + \frac{|\delta^+ + \delta^-|}{|\delta^+|+|\delta^-|}\right)
      \| T_m\otimes S: {L^p_X(\real^d)\to L^p_Y(\real^d)} \|
\end{eqnarray*}
where the equality follows from Lemma \ref{le:multiplier-torus-rn}.
\halmos
\end{proof}
  
\prop{pr:3.2} 
For  $T_{m_0}=-{\rm Id}- 2 R_1^2$, $p\in (1,\infty )$, and $d\geq 2$ one has
\beqla{eq:uusi}   
        \| T_{m_0} \otimes S : L^p_X(\real^d) \rightarrow L^p_Y(\real^d) \|
   \le  \big \| S | (\bsymm{d} ,W) \big \|_p.
\eeq
\eprop
\begin{proof} 
We  apply the  representation of products of Riesz-transforms
in terms of heat extensions to the upper half space (see Lemma \ref{le:A.2}).
Let $f=\sum_{k=1}^m f_k x_k$ and $g=\sum_{l=1}^n g_l b_l$ with 
$f_k,g_l\in \smooth_0 (\real^d)$ and $x_k\in X$, $b_l\in Y'$.
Assume that $u_k$ and $v_l$ are the heat extensions of $f_k$ and $g_l$,
respectively, to the upper half plane and that $(W_t)_{t\ge 0}$ is a standard Brownian motion in  $\real^d$ starting in the origin. 
Let $u := \sum_{k=1}^m u_k x_k$, 
$v:= \sum_{l=1}^n v_l b_l$, and $1=(1/p) + (1/p')$.
Lemma \ref{le:A.2} gives that 
\begin{eqnarray*}
&   & \hspace*{-2em}
      \left|\int_{\real^d} \langle ((T_{m_0} \otimes S)f)(x),g(x)\rangle\, 
      dx\right| \\
& = & \lim_{T\to \infty}(2\pi T)^{d/2}  \bigg | \sum_{k,l} \langle Sx_k,b_l \rangle 
      \times \\
&   & \times \expec \left(\int_{0}^T \nabla  u_k(W_t,T-t)\cdot d (\bsymm d W)_t
                   \int_{0}^T \nabla  v_l(W_t,T-t)\cdot dW_t\right) \bigg | \\
& = & \lim_{T\to \infty}(2\pi T)^{d/2}  \\
&   &  \bigg | \expec \left \langle 
            \int_{0}^T \nabla  S u(W_t,T-t)\cdot d (\bsymm d W)_t,
            \int_{0}^T \nabla    v(W_t,T-t)\cdot dW_t
            \right \rangle \bigg |.
\end{eqnarray*}
We continue with 
\begin{eqnarray*}
&    & \bigg | \expec \left \langle 
            \int_{0}^T \nabla  S u(W_t,T-t)\cdot d (\bsymm d W)_t,
            \int_{0}^T \nabla    v(W_t,T-t)\cdot dW_t
            \right \rangle \bigg | \\
&\leq& \left ( \expec \left \| \int_{0}^T 
       \nabla  S u (W_t,T-t) \cdot d(\bsymm d W)_t \right \|^p_{Y} \right)^{1/p} \\
&    & \hspace*{9em} \times 
       \left( \expec \left \|  \int_{0}^T \nabla  v(W_t,T-t) \cdot dW_t  
       \right \|^{p'}_{Y'} \right)^{1/p'}\\
&\leq& \ito{S}{\bsymm{d}}{p} 
       \left(\expec \left \| \int_{0}^T \nabla  
       u (W_t,T-t) \cdot dW_t \right \|^p_{Y} \right)^{1/p} \\
&    & \hspace*{9em} \times
       \left(\expec \left \| 
       \int_{0}^T \nabla  v(W_t,T-t) \cdot dW_t 
       \right \|^{p'}_{Y'} \right)^{1/p'} \\
& =  & \ito{S}{\bsymm{d}}{p} 
       \| f(W_T) - u(0,T) \|_{L^p_Y} 
       \| g(W_T) - v(0,T) \|_{L^{p'}_{Y'}}
\end{eqnarray*}
by It\^o's formula because 
$(1/2) \Delta u_k = (\partial/\partial t) u_k$ and 
$(1/2) \Delta v_l = (\partial/\partial t) v_l$.
Next, $\sup_{T\ge 0} \| T^{d/2}u(0,T)\|_{X}<\infty$ gives
$\lim_{T\to\infty} T^{d/2} \|  u(0,T) \|^p_X = 0$
so that
$\lim_{T\to\infty} (2\pi T)^{d/2} \expec \| f(W_T) - u(0,T) \|^p_X
 =  \| f\|_{L^p_X(\real^d)}^p$. The same applies for $g(W_T)$ and we end up with
\[  \left|\int_{\real^d} \langle ((T_{m_0} \otimes S)f)(x),g(x)\rangle\, 
     dx\right| 
    \le   \ito{S}{\bsymm{d}}{p} \| f\|_{L^p_X(\real^d)} \| g\|_{L^{p'}_{Y'}(\real^d)}. \]
The proof is complete because $f$ and $g$ as above are dense in 
$L^p_X(\real^d)$ and $L^{p'}_{Y'}(\real^d)$, respectively.
\hspace*{0em}\hfill\halmos
\end{proof}
\medskip 

In order to exploit the quantities  $\ito{S}{A}{p}$ in rigorous arguments, one 
needs at some places approximations of 
stochastic integrals by discrete martingales. So for the general
reader's convenience we switch to the simple discretized version 
$\dito{S}{A}{p}$ introduced below.
Just to prove part of our main result one could proceed more directly 
(see Remark \ref{rem:direct}).   
Given $p\in (1,\infty)$, $S\in L(X,Y)$, and $A\in M(d,d)$, we let
\[
      \dito{S}{A}{p} := \inf c,
\]
where the infimum is taken over all $c>0$ such that 
\begin{multline*}
      \Big \| \sum_{k=1}^N \sum_{l=1}^d 
        \left [ Sd_{(k-1)d+l}(\varphi_1,...,\varphi_{(k-1)d}) \right ]
        \big \langle A (\varphi_{(k-1)d+1},...,\varphi_{kd}), e_l \big \rangle 
                               \Big \|_{L^p_Y} \\
      \le c \Big \| \sum_{k=1}^N   \sum_{l=1}^d
          \left [ d_{(k-1)d+l}(\varphi_1,...,\varphi_{(k-1)d}) \right ]
          \varphi_{(k-1)d+l} 
          \Big \|_{L^p_X},
\end{multline*}
where $N=1,2,...$, 
$d_j:\real^{\left \lfloor \frac{j-1}{d}\right\rfloor d}\to X$ 
are continuous bounded functions taking values in a finite dimensional 
subspace of $X$, and $\varphi_1,\varphi_2,...$ are independent standard 
Gaussian random variables. The $Nd$ terms on the right-hand side 
(in their natural order) we call a Gaussian block martingale difference 
sequence of order $d\ge 1$, on the left-hand side we have its  
{\it $A$-transform}. Observe that
these transforms are not the traditional martingale transforms appearing in
the definition of the UMD-spaces. 

\lem{le:properties_Gaussian_martingales}
Let $A=[a_{lk}]\in M(d,d)$,  
$S\in L(X,Y)$, and $p\in (1,\infty)$.
Then the following is true:
\begin{enumerate}[{\rm (i)}]
\item $\dito{S}{A}{p}=\ito{S}{A}{p}.$
\item If $U\in M(d,d)$ is real and unitary, then 
      \[ \dito{S}{U^TAU}{p}=\dito{S}{A}{p}. \]
\item If $M\geq 1$ is an integer and the tensor product
      $\otimes^M A$ is defined as the block diagonal matrix with $A$
      as each diagonal block, then
      \[ \dito{S}{\otimes^M A}{p}=\dito{S}{A}{p}. \]
\item Assume a sub-matrix $B$ of $A$ obtained from $A$ by choosing 
      indices $I=\{ k_1,\ldots,k_{d'}\}$ with
      $1\leq k_1<k_2<\ldots < k_{d'} \leq d$ and $1\le d' < d$,
      and deleting the corresponding rows and columns from $A$.
      Then
      \[ \dito{S}{B}{p}\leq\dito{S}{A}{p}. \]
\end{enumerate}
\elem
\begin{proof}
(i) It is evident that $\dito{S}{A}{p}\leq\ito{S}{A}{p}.$ The
inequality to the other direction follows by a standard approximation
of Ito-integrals by discrete Gaussian martingales.
\medskip

(ii) Here we observe that 
\begin{eqnarray*}
&   & \hspace*{-1em}
      \Big \| \sum_{k=1}^N \sum_{l=1}^d
      \left [ Sd_{(k-1)d+l}(\varphi_1,...,\varphi_{(k-1)d}) \right ]
      \big \langle U^TAU(\varphi_{(k-1)d+1},...,\varphi_{kd}), e_l \big \rangle 
      \Big \|_{L^p_Y} \\
& = & \Big \| \sum_{k=1}^N 
      \Big \langle U \big ( (Sd_{(k-1)d+l}(\varphi_1,...,\varphi_{(k-1)d}))_{l=1}^d
                     \big ), AU (\varphi_{(k-1)d+1},...,\varphi_{kd})
       \Big \rangle \Big \|_{L^p_Y} \\
&\le& \dito{S}{A}{p} \;\; \cdot\\
&   & \Big \| \sum_{k=1}^N 
       \Big \langle U \big ( (d_{(k-1)d+l}(\varphi_1,...,\varphi_{(k-1)d}))_{l=1}^d
                     \big ),
                    U(\varphi_{(k-1)d+1},...,\varphi_{kd})
       \Big \rangle \Big \|_{L^p_Y} \\
& = & \dito{S}{A}{p} 
        \Big \|  \sum_{k=1}^N \sum_{l=1}^d 
          \left [ d_{(k-1)d+l}(\varphi_1,...,\varphi_{(k-1)d}) \right ]
          \varphi_{(k-1)d+l} \hspace{.5em}
          \Big \|_{L^p_X}.
\end{eqnarray*}
Above we used in the second step the observation that
the transformed sequence
\beqla{eq:3.400}
&   & \langle U(\varphi_{1},...,\varphi_{d}),e_1\rangle,
      \ldots,
      \langle U(\varphi_{1},...,\varphi_{d}),e_d\rangle, \nonumber\\
&   & \langle U(\varphi_{d+1},...,\varphi_{2d}),e_1\rangle,
      \ldots,
      \langle U(\varphi_{d+1},...,\varphi_{2d}),e_d\rangle, \ldots \nonumber
\eeq
consists again of independent Gaussian random variables because $U$ is unitary. 
\medskip

(iii) The $\otimes^M A$-transform of an appropriate
Gaussian martingale block difference sequence of order $dM$ 
is obtained by simply performing the $A$-transform
of the same sequence (which is also of order $d$). 
This shows that $\dito{S}{\otimes^M A}{p}\leq \dito{S}{A}{p}.$
The converse inequality is a special case of (iv) we treat now.

\medskip

(iv) We first perform a unitary permutation of
the coordinates, which is justified by part (ii),
so that we may assume that $B=(a_{lk})_{1\leq l,k\leq d'}$. 
Next we consider the unitary map $U$ which maps the unit vector $e_k$
to $e_k$ if $k=1,...,d'$ and to $-e_k$ if $k=d'+1,...,d$. From (ii) it
follows that, for $C:= (1/2) (A + U^T A U)$, one has
$\dito{S}{C}{p}\leq \dito{S}{A}{p}$. The entries of $C$ satisfy
$c_{lk}=a_{lk}$ if $l,k \in  \{ 1,...,d'\}$ and $c_{lk}=0$ if $l\in \{1,...,d'\}$
and $k\not\in \{ 1,...,d'\}$.
Now the inequality 
$\dito{S}{B}{p}\leq \dito{S}{C}{p}$
can be proved by an appropriate augmentation of the
Gaussian random variables and of the martingale difference sequence
that has to be transformed: we add to each block 
$\phi_{(k-1)d'+1},\ldots,\phi_{kd'}$ of Gaussian random variables  
$d-d'$ independent Gaussian random variables to obtain a block size $d$
and add appropriate zero martingale differences to the original 
martingale difference sequence. Now the transformation with respect to $B$ 
can be artificially written as a transformation with respect to $C$.
\halmos
\end{proof}
\bigskip

The previous lemma will be used to prove the following result:

\prop{pr:3.3}
{\rm (i)}\quad Assume that $A\in M(d,d)$ is real and symmetric, and 
denote by
 $\lambda_{\rm max}$ (respectively, $\lambda_{\rm min}$) the largest
(respectively, the smallest) eigenvalue of $A$. Let $B\in M(d',d')$ be
another real and symmetric matrix such that each eigenvalue $\lambda$
of $B$ satisfies $\lambda_{\rm min}\leq\lambda\leq\lambda_{\rm max}.$
Then for any $S\in L(X,Y)$ and $p\in (1,\infty ),$ it holds that 
$$
\ito{S}{B}{p}\leq\ito{S}{A}{p}.
$$
\noindent{\rm (ii)}\quad Assume that $A\in M(d,d)$ is real and 
antisymmetric, 
and denote by $\lambda_{\rm max}$ (respectively, $\lambda_{\rm min}$) the largest
(respectively, the smallest) eigenvalue of $iA$. Let $B\in M(d',d')$ be
another real and antisymmetric matrix such that each eigenvalue $\lambda$
of $iB$ satisfies $\lambda_{\rm min}\leq\lambda\leq\lambda_{\rm max}.$
Then for any $S\in L(X,Y)$ and $p\in (1,\infty ),$ it holds that 
$$
\ito{S}{B}{p}\leq\ito{S}{A}{p}.
$$
\eprop
\begin{proof}
(i) \quad By Lemma \ref{le:properties_Gaussian_martingales} (i) it is enough to show that
$\dito{S}{B}{p}\leq\dito{S}{A}{p}$. Moreover, by applying the
spectral theorem for symmetric matrices and part (ii) of the same
Lemma we may assume that both $A$ and $B$ are (real and) diagonal.
In addition, by Lemma \ref{le:properties_Gaussian_martingales} (iii) we may replace $A$
by the tensor product $\otimes^{d'} A$. Observe that the tensor product
$\otimes^{d'} A$ is diagonal and at least $d'$ of its diagonal
elements have the value $\lambda_{\rm max}$, and the same 
holds true for $\lambda_{\rm min}$. By applying again a unitary
permutation for the coordinates and part (iv) of Lemma \ref{le:properties_Gaussian_martingales}
we obtain for any given sequence $\Lambda := (\lambda_1,\ldots ,\lambda_{d'})$
with $\lambda_j\in\{ \lambda_{\rm max}, \lambda_{\rm min}\}$ for all $1\leq j\leq d'$
that the diagonal matrix  $A_\Lambda\in M(d',d')$ (with the diagonal $\Lambda$)
satisfies
$$
\ito{S}{A_\Lambda}{p}\leq\ito{S}{A}{p}.
$$ 
By the assumption on the eigenvalues of
the matrix $B$, we may express $B$ as a convex combination of
matrices of the form $A_\Lambda$. This clearly yields the claim.

\medskip 

(ii)\quad The matrix $iA$ is self-adjoint, so that the  eigenvalues are real.
A simple examination of the spectral decomposition of $iA$ (observe that
$\lambda$ and $-\lambda$ are simultaneously eigenvalues for $iA$)
we may write $A$, after a unitary transformation, in the form
$$
A= B_1\otimes\ldots \otimes B_\ell,
$$
where each $B_k$ is of the form $B_k=c_k\basym$ 
($\basym$ was defined in (\ref{eqn:matrices})), with 
$\lambda_{\min} \le c_k \le \lambda_{\max}$
(in case $d$ is odd the last one, i.e. $B_\ell$,
equals the $1\times 1$ zero matrix and 
$\lambda_{\min} \le 0 \le \lambda_{\max}$).
Now the claim follows by an extreme point argument like in (i).
\halmos
\end{proof}
\bigskip
 
\prop{pr:Tm-Tm0}
Let $m$ be a homogeneous,  even, and smooth multiplier on $\real^d,$ $d\geq 2.$ Then 
for any $S\in L(X,Y)$ it holds that
\[      \left \| T_m\otimes S:L^p_X(\real^d)\to L^p_Y(\real^d)\right \|
   \leq c \hspace{.1em}
        \left \| T_{m_0}\otimes S:L^p_X(\real^d)\to L^p_Y(\real^d) \right \| \]
for  $S\in L(X,Y)$ and $p\in (1,\infty)$, where $c>0$ depends 
at most on $m$.
\eprop
\begin{proof} 
In the following we always assume that $\xi,\theta\in S^{d-1}$, so that
(for example) $m_0(\xi )= 2\xi_1^2-1$. Let $a\in (0,1)$. 
By composing  $m_0$ with the linear map 
$\xi\mapsto B_a(\xi):=(\sqrt{1-a^2}\xi_1,\xi_2,\ldots, \xi_d)$ we infer 
by \refeq{eq:2.10} that 
\[     \|T_{m_a}\otimes S: L_X^p(\real^d) \to L_Y^p(\real^d) \|
   \le \frac{1}{1+a} 
       \|T_{m_0}\otimes S: L_X^p(\real^d) \to L_Y^p(\real^d) \|, \]
where
\begin{eqnarray*}
      m_a(\xi )
&:= & {2-a^2\over 2(1+a)}+{a^2\over 2(1+a)}m_0\circ B_a(\xi )=(1-a){1\over 1-a^2\xi_1^2}
      \nonumber\\
& = &(1-a)\sum_{k=0}^\infty a^{2k}\xi_1^{2k}
\end{eqnarray*}
and we also use the estimate 
\[       \| Id\otimes S : L_X^p(\real^d) \to L_Y^p(\real^d) \|
    \leq \| S \| 
    \leq \|T_{m_0}\otimes S: L_X^p(\real^d) \to L_Y^p(\real^d) \| \]
(for example, use Lemma \ref{le:multiplier-torus-rn} and $\| m_0 \|_\infty =1$).
Of course, given $\theta\in S^{d-1}$ we obtain the same estimates for 
\[     m_a^\theta(\xi )
   := (1-a)\sum_{k=0}^\infty a^{2k}\langle \xi,\theta \rangle^{2k}. \]
Given $r\in (0,1)$, the $d$-dimensional Poisson kernel for $B(0,1)$ has the form
\beqla{eq:3.430}
      P(r \xi,\theta)
& = & {1\over \omega_{d-1}}{(1-r^2)\over (1+r^2-2r\langle \xi,\theta\rangle)^{d/2}}
      \nonumber\\
& = & {(1-r^2)(1+r^2)^{-d/2}\over \omega_{d-1}} 
      \left (1-{2r\langle \xi,\theta \rangle\over 1+r^2}\right )^{-d/2} \nonumber
\eeq
with $\omega_{d-1}= |\partial B(0,1)|$.
By the Taylor expansion of $y\mapsto (1-y)^{-d/2}$ we obtain that
\beqla{eq:3.435}
      P^s(r\xi,\theta)
&:= & {1\over 2}(P(r\xi,\theta)+P(r\xi,-\theta))\nonumber\\
& = & {\frac{(1-r^2)(1+r^2)^{-d/2}}{\omega_{d-1}}}
      \sum_{k=0}^\infty \left({\frac{2r}{1+r^2}}\right)^{2k}
      \langle \xi,\theta \rangle^{2k}{\binom{-d/2}{2k}}. \nonumber
\eeq
Fix $\varepsilon >0$ and let 
\[ c_0 :=  2  \omega_{d-1} \frac{\Gamma(\varepsilon+1) \Gamma(d/2)}
                                {\Gamma((d/2)+\varepsilon)}. \]
Using the substitution $u=2r/(1+r^2)$ and Euler's $\beta$-integral yields that
\beqla{eq:3.440}
&   & \int_0^1(1+r^2)^{-2+d/2}
      \left (1-{2r\over 1+r^2}\right )^{d/2-1+\varepsilon}P^s(r\xi,\theta)\, dr
      \nonumber\\
& = & \frac{1}{2\omega_{d-1}} 
      \sum_{k=0}^\infty \left (\int_0^1(1-u)^{d/2-1+\varepsilon}u^{2k}\,du 
      \right ){\binom{-d/2}{2k}}\langle \xi,\theta \rangle^{2k} \nonumber\\
& = & \frac{1}{2\omega_{d-1}} {\frac{\Gamma (d/2 + \varepsilon )}{\Gamma (d/2)}}\sum_{k=0}^\infty 
      \left ({\Gamma (d/2 + 2k )\over\Gamma (d/2+2k+1+\varepsilon )} \right )
      \langle \xi,\theta \rangle^{2k} \nonumber \\
& = & \frac{1}{c_0} 
      \int_0^1(1-a)^{\varepsilon -1}a^{d/2-1}m_a^\theta(\xi)\, da.
\eeq
Let us denote by $\sphe(d)$ a complete orthonormal
$L^2$-basis of $L^2(S^{d-1}, \lambda)$, where $\lambda$ is the normalized
Haar measure, consisting of spherical 
harmonics on $S^d$ (we refer to \cite{Stein2,SteinWeiss} for more details 
on spherical harmonics). For $\psi\in \sphe(d)$ we denote by
${\rm deg}(\psi)$ the degree of $\psi$. Because the function 
$(r,\xi)\mapsto r^k\psi (\xi)$ is harmonic if ${\rm deg}(\psi )=k$
we have in  this case
\begin{equation}\label{eq:eigenvalue_P}
     \int_{S^{d-1}}P(r\xi ,\theta) 
     \psi(\theta)\, d\theta
   = r^k \psi(\xi ). 
\end{equation}
Define
\[    \lambda_k 
   := c_0 \int_0^1 (1+r^2)^{-2+d/2}
          \left (1-{2r\over 1+r^2}\right )^{d/2-1+\varepsilon}r^k\,dr 
  \sim_{c(\varepsilon ,d)}k^{1-d-2\varepsilon}. \]
Let $a_\psi:= \langle m,\psi \rangle/\lambda_k$ if $\psi\in \sphe(d)$ and
${\rm deg}(\psi )=k$. Applying \cite[p. 70]{Stein2} gives 
\[ \sum_{\psi\in \sphe_d} |a_\psi|^2 < \infty. \]
Because $m$ is even we find an even $f\in L^2(S^{d-1},\lambda)$ with
$\langle f,\psi \rangle = a_\psi$, so that
\begin{equation}\label{eq:representation_m}
m = \sum_{k=0}^\infty \lambda_k \sum_{\psi\in \sphe_d,\; {\rm deg}(\psi )=k}
       \langle f,\psi\rangle \psi
\end{equation}
in $L^2(S^{d-1},\lambda)$. If we can show that
\begin{equation}\label{eq:integral_representation_m}
     m(\xi) 
   = \int_{S^{d-1}} \int_0^1 (1-a)^{\varepsilon-1}a^{\frac{d}{2}-1}
     m_a^\theta(\xi) f(\theta) da d\theta
\end{equation}
for $\xi\in S^{d-1}$, then this would give our assertion since
\[ \int_{S^{d-1}}\int_0^1  (1-a)^{\varepsilon-1} a^{\frac{d}{2}-1} |f(\theta)| da d\theta
   < \infty \]
implies that $m$ is a convex combination of the multipliers $m_a^\theta$.
In order to verify equality (\ref{eq:integral_representation_m})
it is sufficient to show that
\begin{multline*}
    \int_{S^{d-1}} \left [ \int_{S^{d-1}} 
    \int_0^1 (1-a)^{\varepsilon-1}a^{\frac{d}{2}-1}
    m_a^{\theta}(\xi) f(\theta) da d\theta \right ] \overline{\psi}(\xi) d \xi \\
  = c_0 \int_0^1 (1+r^2)^{-2+d/2}
    \left (1-{2r\over 1+r^2}\right )^{d/2-1+\varepsilon}r^k\,dr 
    \langle f,\psi\rangle
\end{multline*}
for $\psi\in \sphe_d$ with ${\rm deg}(\psi )=k$.
But this follows by a computation from (\ref{eq:3.440}) and 
(\ref{eq:eigenvalue_P}).
\halmos
\end{proof}

\begin{remark}
A slight additional argument shows that one may allow any even multiplier $m$ in the above 
proposition such that $m_{|S^{d-1}}\in W^{s ,1}(S^{d-1})$ for some $s> d-1.$ Especially, if
$d=2$ this class contains all $m$ such that for some $\varepsilon >0$ the derivative 
$D^\varepsilon m_{|S^{1}}$ is a function of bounded variation.
\end{remark}
\medskip
We are ready for
\bigskip

\noindent{\bf Proof of Theorem \ref{th:3.1}.}\quad
We first verify that $\umd_p(S)$ is comparable to the corresponding norm of 
the operator $T_{m_0}\otimes S.$

Observe  that $\max_{|\xi | =1} m_0(\xi )=1$ and $\min_{|\xi | =1} m_0(\xi )=-1$.
An application of Proposition \ref{pr:UMD-multiplier} yields that
\beqla{eq:3.500}
{\umd}_p(S) 
   \leq  \| T_{m_0}\otimes S  
   : L^p_X(\real^d)  \to L^p_Y(\real^d)\|. 
\eeq 
By Proposition \ref{pr:3.2} we have in turn
\beqla{eq:3.510}
 \| T_{m_0}\otimes S  
   : L^p_X(\real^d)  \to L^p_Y(\real^d)\|\leq \ito{S}{\bsymm{d}}{p}. 
\eeq 
The eigenvalues of $\bsymm{d}$ and $\bsym$ are $\pm 1$. Hence Proposition
\ref{pr:3.3} (i) and Lemma \ref{le:properties_Gaussian_martingales} (i) yield that
\beqla{eq:3.520}
\ito{S}{\bsymm{d}}{p}\leq \ito{S}{\bsym}{p} =\dito{S}{\bsym}{p}.
\eeq
However, the $\bsym$-transform of discrete Gaussian martingales
is a special case of a UMD-martingale transform and
hence in case $S={\rm Id}$ it follows immediately that
\beqla{eq:3.530}
\dito{S}{\bsym}{p}\leq \umd_p (S).
\eeq
There are various ways to check this inequality for general $S$ (note that we started with
Walsh-Paley martingales in the definition of $\umd_p (S)$). 
An easy self-consistent way would be to apply
the central limit theorem argument from Lemma \ref{lemma:lower_bound_HT} to the 
Bernoulli variables to replace the Bernoulli variables by the Gaussian random variables
(there are arguments to switch from Walsh-Paley martingale difference sequences to 
arbitrary martingale difference sequences; see for example \cite[p. 12]{Burk4}
and \cite{Maurey}).
By combining the inequalities \refeq{eq:3.500}--\refeq{eq:3.530}
we obtain  that 
\beqla{eq:3.540}
&&   \umd_p (S) 
   = \| T_{m_0}\otimes S : L^p_X(\real^d)  \to L^p_Y(\real^d)\|
   = \ito{S}{\bsym}{p}.
\eeq
This proves the equalities in the theorem.

Now let us assume that $A\in M(d,d)$ is symmetric and non-trivial, 
and that $m$ is an even, smooth, and non-trivial multiplier. 
Proposition \ref{pr:UMD-multiplier} shows that $\umd_p(S)$ is
dominated by a multiple of the norm of $T_m\otimes S$, and 
Proposition \ref{pr:Tm-Tm0} verifies that this norm in turn is dominated by 
a multiple of the norm of $T_{m_0}\otimes S$. The linear 
equivalence of the norm of $T_m\otimes S$
to the UMD-constant of $X$ follows now from \refeq{eq:3.540}.  
 
Finally, the equivalence of $\umd_p(S)$ to $\ito{S}{A}{p}$  will be deduced 
from \refeq{eq:3.540} by an application of 
Proposition \ref{pr:3.3} (i). Denote by $\lambda_{\rm min} , \lambda_{\rm max}$ 
the minimal and maximal eigenvalue of $A$.
Firstly, we obtain that 
\[ \ito{S}{A}{p}\leq c\ito{S}{\bsym}{p} \]
with 
 $c:=\max (|\lambda_{\rm max}|,|\lambda_{\rm min}| )$.
To consider the other direction, let 
$\alpha:=2/(\lambda_{\rm max}-\lambda_{\rm min})$
and 
$\beta:=-(\lambda_{\rm max}+\lambda_{\rm min})/
         (\lambda_{\rm max}-\lambda_{\rm min})$, so that
$B:= \alpha A +\beta$ satisfies  
$\lambda_{\rm max}(B)=-\lambda_{\rm min}(B)=1$.
Consequently,
\begin{eqnarray*}
      \ito{S}{\bsym}{p}
&\le& \alpha \ito{S}{A}{p} + \beta \ito{S}{Id}{p} \\
&\le& \alpha \ito{S}{A}{p} 
       + \frac{\beta}{\max\{ |\lambda_{\max}|,|\lambda_{\min}|\} }
         \ito{S}{A}{p}
\end{eqnarray*}
and the proof of Theorem \ref{th:3.1} is complete.
\hfill\halmos
\bigskip 

\begin{remark}\label{rem:direct}
There is a shorter argument for part of Theorem \ref{th:3.1} which, 
modulo some 'hand waving', bypasses Lemma
\ref{le:properties_Gaussian_martingales} and Proposition \ref{pr:3.3}. 
Namely, in order to prove e.g. the statement
\[
{\rm UMD}_p(X)
    = \| T_{m_0}  : L^p_X(\real^d)\to L^p_Y(\real^d) \|
   \sim_C \| T_m : L^p_X(\real^d)\to L^p_Y(\real^d) \|
\]
one only needs Propositions \ref{pr:UMD-multiplier}, \ref{pr:3.2}, and \ref{pr:Tm-Tm0} and 
$\ito{S}{\bsymm{d}}{p} \le \umd_p(S)$. Because $\bsymm{d}$ is a
diagonal matrix, the latter inequality can be suitably approximated
by a UMD-transform inequality. 
\end{remark}

%%%%%%%%%%%%%%%%%%%%%%%%%%%%%%%%%%%%%%%%%%%%%%%%%%%%%%%%%%%%%%%%%%%%%%%%%%%%%

\section{Odd multipliers}\label{se:hilbert}

A canonical example of an odd multiplier is the Hilbert transform. In this section we
verify that all odd multipliers, as well as all antisymmetric transforms for stochastic
integrals are linearly comparable to the Hilbert transform. Observe that an odd
multiplier maps real functions to purely imaginary ones, whence we  
consider below the operator $T_{im}$ in order to be able 
to allow also real Banach spaces
in the statement.
 
\thm{th:4.1}
Assume that $m\in\smooth (\real^d\backslash \{0\} )$, $d\ge 2$, is a real non-zero, 
homogeneous, and odd multiplier, and that $A\in M(d,d)$ is a non-zero and 
anti-symmetric matrix. Let $p\in (1,\infty ).$ Then there is a constant  $C=C(A,m)$  
such that for every pair of Banach spaces $X$ and $Y$ and for every operator 
$S\in L(X,Y)$ one has 
\begin{eqnarray*}
         \| T_{im}\otimes S : L^p_X(\real^d)\to L^p_Y(\real^d) \|  
&\sim_C& \ito{S}{A}{p} \sim_C \ito{S}{\basym}{p} \\
&  =   & \| \htt\otimes S : L^p_X(\real )\to L^p_Y(\real ) \|.
\end{eqnarray*}
\ethm

We start the proof by
\lem{lemma:lower_bound_HT}
For $p\in (1,\infty)$ and $S\in L(X,Y)$ one has that 
\[     \ito{S}{\basym}{p}
   \leq\| 
   \widetilde\htt\otimes S : L^p_X(\torus)\to L^p_Y(\torus) \|. \]
\elem

\begin{proof}
Assume that  $M\geq 1$ and that $a_k,b_k: [-1,1]^{2k-2} \rightarrow X$,
$k=1,...,M$, are bounded and continuous. Lemma \ref{le:Bourgain_new} yields 
that
\beqla{eq:4.20}
&   & \Big \| \sum_{k=1}^{M}
      \Big ( \widetilde\htt \sin(\theta_k) S a_k(\sin(\theta_1),...,\sin(\theta_{k-1});
                                                 \cos(\theta_1),...,\cos(\theta_{k-1})) \\
&   & + \widetilde\htt\cos(\theta_k) S 
             b_k(\sin(\theta_1),...,\sin(\theta_{k-1});
                 \cos(\theta_1),...,\cos(\theta_{k-1}) \Big ) \big \|_{L_p^Y}\nonumber \\
&\le& \|\widetilde{\htt}\otimes S:L_p^X(\torus) \rightarrow L_p^Y(\torus)\|\; 
      \times \nonumber\\
&   & \times \Big \| \sum_{k=1}^{M}
      \Big ( \sin(\theta_k) a_k(\sin(\theta_1),...,\sin(\theta_{k-1});
                           \cos(\theta_1),...,\cos(\theta_{k-1})) \nonumber\\
&   & + \cos(\theta_k) b_k(\sin(\theta_1),...,\sin(\theta_{k-1});
                           \cos(\theta_1),...,\cos(\theta_{k-1}) \Big ) 
      \big \|_{L_p^X}.\nonumber
\eeq
Next, we apply a blocking argument. Let $N, L\geq 1$ be integers and assume that
$A_k,B_k: [-1,1]^{2k-2} \rightarrow X$, $k=1,...,N$, are bounded and continuous
functions. For $M=NL$ we apply \refeq{eq:4.20} to 
\beqla{eq_4.30}
&   & a_{(k-1)L+j}(\sin(\theta_1),...,\sin(\theta_{(k-1)L+j-1});
                           \cos(\theta_1),...,\cos(\theta_{(k-1)L+j-1}))\nonumber\\
& = & A_k\bigg( (L/2)^{-1/2}(\sin(\theta_1)     +\ldots + \sin (\theta_L)), \nonumber \\
&   & \hspace*{4em} (L/2)^{-1/2} (\sin(\theta_{L+1}) +\ldots + \sin (\theta_{2L})),\nonumber\\
&   & \hskip1truein\ldots ,(L/2)^{-1/2}
                       (\sin (\theta_{(k-2)L+1})+\ldots + \sin (\theta_{(k-1)L})); \nonumber\\
&   & (L/2)^{-1/2}(\cos(\theta_1)+\ldots + \cos (\theta_L)), \\
&   & \hspace*{4em} (L/2)^{-1/2} (\cos (\theta_{L+1})+\ldots + \cos (\theta_{2L})), \nonumber\\
&   & \hskip1truein \ldots ,(L/2)^{-1/2}(\cos(\theta_{(k-2)L+1})+\ldots 
      + \cos (\theta_{(k-1)L}))\bigg) \nonumber
\eeq
and with an analogous choice for the coefficients $b_{(k-1)L+j}$ for $1\leq k\leq N$ 
and $1\leq j\leq L$. By taking into account the action of $\widetilde{\htt}$ on the
trigonometric polynomials we obtain that
\beqla{eq:4.50}
&   & \bigg \| \sum_{k=1}^{N}
      \Big ( -S c_k A_k(s_1,...,s_{k-1};
                           c_1,...,c_{k-1})\; +\nonumber\\
&   & \hskip1truein +\; S s_k B_k(s_1,...,s_{k-1};
                           c_1,...,c_{k-1}) \Big ) \bigg \|_{L_p^Y}\nonumber \\
&\le& \|\widetilde{\htt}\otimes S:L_p^X(\torus) \rightarrow L_p^Y(\torus)\|\;\times 
      \nonumber\\
&   & \times \bigg \| \sum_{k=1}^{N} \Big ( s_k A_k(s_1,...,s_{k-1};
                           c_1,...,c_{k-1}) \; +\nonumber\\
&   & \hskip1truein +\; c_k B_k(s_1,...,s_{k-1};
                           c_1,...,c_{k-1}) \Big ) \bigg \|_{L_p^X}\nonumber
\eeq
where we have used the abbreviations 
$s_k:=(L/2)^{-1/2}\sum_{j=1}^L\sin (t_{(k-1)L+j})$ and 
$c_k:=(L/2)^{-1/2}\sum_{j=1}^L\cos (t_{(k-1)L+j})$ for $k=1,\ldots ,N.$
By keeping $N$ fixed, letting $L\to\infty$, and by applying the central
limit theorem (here we may normalize the Lebesgue measure on $\torus$)
and the fact that $\cos$ and $\sin$ are uncorrelated, one gets
that 
\[      \dito{S}{\basym}{p} 
   \le \| \widetilde\htt\otimes S:L^p_X(\torus) \rightarrow L^p_Y(\torus) \|.  \]
Finally, Lemma \ref{le:properties_Gaussian_martingales} verifies that 
$\dito{S}{\basym}{p} =\ito{S}{\basym}{p}$ and
we are done.
\hspace*{0em}\hfill\halmos
\end{proof} 

The following lemma is well-known. 
For the convenience of the reader we recall the idea of its proof.

\lem{lemma:upper_bound_HT}
For $p\in (1,\infty)$ one has that 
\[
       \| \htt\otimes S : L^p_X(\real )\to L^p_Y(\real ) \| 
    \le  \ito{S}{\basym}{p}.
\]
\end{lemma}
\begin{proof}
We use
$  \| \htt\otimes S             : L^p_X(\real )    \to L^p_Y(\real) \|  
 = \| \widetilde{\htt}\otimes S : L^p_{X,0}(\torus)\to L^p_Y(\torus)\|$
from Lemma \ref{le:multiplier-torus-rn} and consider
$   f(\theta) 
 := \sum_{k=1}^n \left ( \sin (k\theta) x_k + \cos (k\theta) y_k \right )$
with $x_k,y_k \in X$ as a function on the unit circle. Then, a.s.,
\[   u(W_{\tau} ) 
   = \int_0^\tau \nabla u(W_t)\cdot dW_t
 \sptext{.5}{and}{.5}
     v(W_\tau) 
   = - \int_0^\tau         \nabla u(W_t)\cdot d(\basym W)_t \]
by It\^o's formula, where $u$ and $v$ are the harmonic extensions 
of $f$ and
$\sum_{k=1}^n \left ( - \cos (k\theta) x_k + \sin (k\theta) y_k \right )$
to the unit disc, $(W_t)_{t\ge 0}$ is a standard two-dimensional standard 
Brownian motion, and 
$\tau := \inf \{ t \ge 0 : |W_t|=1 \}$. 
\hfill\halmos
\end{proof} 
\bigskip 

{\bf Proof of Theorem \ref{th:4.1}}.\quad 
The equality
\[  \ito{S}{\basym}{p} 
  = \| \htt\otimes S : L^p_X(\real )\to L^p_Y(\real )\| \]
follows from Lemmas \ref{lemma:lower_bound_HT}, \ref{lemma:upper_bound_HT}, and \ref{le:multiplier-torus-rn}.
Moreover, 
\[ \ito{S}{A}{p} \sim_C \ito{S}{\basym}{p} \]
is a consequence of Proposition \ref{pr:3.3} (ii) (note, that the eigenvalues of
$i \basym$ are $\pm 1$ and that $iA$ has at least two symmetric non-zero 
real eigenvalues). Moreover, by a 
classical  argument called 'the method of rotations' (see 
e.g. \cite[p.271, formula (4.2.20)]{Grafakos} one may express any odd (and smooth)
multiplier operator $T_{im}$ as an average of directional 
Hilbert-transforms which immediately yields
that the norm of $T_{im}\otimes S$ is 
linearly dominated by that of $\htt\otimes S$.
Finally, if $m$ is an odd and non-zero multiplier, we may assume that $m(e_1)=1.$
The corresponding discrete multiplier $\widetilde{m}$ satisfies
$-i \widetilde{m}(k_1 e_1) = -i \sgn (k_1)$ for $k_1\not = 0$ and
$-i \widetilde{m}(0 e_1) = 0$. By the consideration of
$    \widetilde{f}(\theta_1,...,\theta_d) 
   = f(\theta_1) 
  := \sum_{k=1}^n 
     \left ( \sin (k_1 \theta_1) x_{k_1} + \cos (k_1\theta_1) y_{k_1} \right )$
with $x_{k_1},y_{k_1} \in X$ and observing that
$  ((\widetilde{H}\otimes S)f)(\theta_1) 
 = ((-T_{i \widetilde{m}}\otimes S)\widetilde{f})(\theta_1,...,\theta_d)$ we 
immediately get that
\[     \| \widetilde{H}\otimes S : L_{X,0}^p(\torus)\to L_Y^p(\torus) \|
   \le \| T_{i\widetilde{m}} \otimes S :  L_X^p(\torus^d)\to L_Y^p(\torus^d) \|\]
so that  
$    \| H      \otimes S : L_X^p(\real)  \to L_Y^p(\real) \|
 \le \| T_{im} \otimes S : L_X^p(\real^d)\to L_Y^p(\real^d) \|$
by  Lemma \ref{le:multiplier-torus-rn}.
\hfill\halmos
\bigskip  

%%%%%%%%%%%%%%%%%%%%%%%%%%%%%%%%%%%%%%%%%%%%%%%%%%%%%%%%%%%%%%%%%%%%%%%%%%%%%

\section{Additional remarks}
\label{se:remarks}

The remaining main open problem is whether or not we have the linear equivalence 
of the norms
\begin{equation}\label{eqn:conjecture}
        \umd_p(S) 
   \sim \| \htt\otimes S : L^p_X(\real)    \to L_Y^p(\real) \|.
\end{equation}
By the results obtained in Theorems \ref{th:3.1} and  \ref{th:4.1}
this problem can be formulated now in 
various ways, for example in a purely probabilistic way via 
martingale transforms or via multipliers. 
To find counterexamples to (\ref{eqn:conjecture})
in the category of Banach spaces (i.e. $S=Id_X$ for a 
Banach space $X$) seems to be 
harder than to find counterexamples for operators $S\in L(X,Y)$,
which is one reason for our usage of the operator setting.
At the same time the operator setting gives an easier control whether
or not estimates are linear 
(for example 
$\umd_p(Id_X)<\infty$ if and only if 
$\| \htt\otimes Id_X : L^p_X(\real)    \to L_X^p(\real) \|<\infty$
because of $Id_X^2=Id_X$ which does not work for operators).
Natural candidates to disprove (\ref{eqn:conjecture}) 
in the setting of operators are the operators of summation 
$\sigma_n:\ell^1_n\to\ell^\infty_n$. For them it is known that 
\[      \| \htt\otimes \sigma_n : L^2_{\ell^1_n}(\real) \to L_{\ell^\infty_n}^2(\real) \|
   \sim \log(n+1), \]
see \cite{DefantM-PHD} (Section IV, Satz 2.1 (proof) and Korollar 2.4)
(cf. also \cite{Pietsch-Wenzel} (2.4.4)).
So far, the best estimates for the $\umd$-constants are
\[     \frac{1}{c}\sqrt{\log(n+1)}
   \le \umd_2(\sigma_n) \le c \log(n+1) \]
where the lower estimate follows from \cite{Geiss8.5} 
and the upper one is a consequence of 
$\umd_2(Id_{\ell^\infty_n}) \le c \log(n+1)$ which is folklore.

A determination of the $\umd$-constant of $\sigma_n$ would be of interest for 
several reasons:
in case of $\umd_2(\sigma_n) \sim  \log(n+1)$ this would imply that each
non-superreflexive Banach space $X$ contains $n$-dimensional subspaces $E_n$
such that the lower bound $\umd_2(Id_{E_n}) \ge \log(n+1)/c$ holds, because
due to R.C. James a Banach space $X$ is non-superreflexive 
if and only if the operators $\sigma_n$ can be uniformly factorized through
certain $n$-dimensional subspaces $E_n\subseteq X$ (see \cite{James}).
Having the lower estimate one might ask for more connections between the
quantitative behavior of the UMD-constants of the finite dimensional subspaces 
of a Banach space and the property that the space is non-super\-reflexive.    
On the other hand, any estimate
of type $\umd_2(\sigma_n) = o(\log(n+1))$ would offer some new inside into
martingale transforms. 

%%%%%%%%%%%%%%%%%%%%%%%%%%%%%%%%%%%%%%%%%%%%%%%%%%%%%%%%%%%%%%%%%%%%%%%%%%%%%%

\section{Appendix}\label{se:appendix}

We recall the well-known connection between 
singular integrals and $A$-transforms of stochastic integrals  corresponding 
to heat extensions of functions to the upper half space. A simple proof will be sketched 
below for the readers convenience. We refer to \cite{Ban2} for 
the original result, and to \cite{Bass} and the references therein for corresponding 
results that use the harmonic extension instead. In order 
to state the relation let $d\geq1$ and assume a real matrix 
$A=[a_{kl}]\in M(d,d)$.
We define the operator $U_A$ for smooth elements 
$f,g\in \smooth_0 (\real^d)$  through the bilinear form
\beqla{eq:A.100}
&   & \int_{\real^d} (U_Af)(x)g(x)\, dx \\
&:= & \lim_{T\to \infty}(2\pi T)^{d/2}\expec \left[\left(\int_{0}^{T} 
      \nabla  u(W_t, T-t)\cdot d(AW)_t \right)\times \right. \nonumber \\
&   & \qquad \times  \left. \left(\int_{0}^{T}
      \nabla  v(W_t, T-t)\cdot dW_t\right)\right] ,\nonumber 
\eeq
where $u$ and $v$ are the heat extensions of $f$ and $g$, respectively,
to the upper half space $\real^{d}\times \real_+ $ (i.e. $u(t,x) := \E f(x+W_t)$ and
similarly for $v$), $\nabla$ denotes the differentiation with respect to the $x$-variables, 
and $(W_t)_{t\ge 0}$ is a standard $d$-dimensional Brownian motion starting at the origin.
Hence, for example, $u_{t}={1\over 2}\Delta u$ and $u(x,0)=f(x)$ for
$x\in\real^d.$  
\bigskip

\begin{lemma}\label{le:A.2}
Let $p\in (1,\infty)$ and $d\geq 2.$ The operator $U_A$ is well-defined and extends 
to a bounded operator on $L^p(\real^d)$ which can be expressed in terms 
of the Riesz transforms as
\[ U_A=-\sum_{k,l=1}^da_{kl}R_kR_l. \]
\end{lemma}
\begin{proof}
We let $f,g\in \smooth_0 (\real^d)$. From It\^o's isometry we obtain that
\begin{eqnarray*}
&   & (2\pi T)^{d/2}\expec \left[\left(\int_{0}^{T} 
      \nabla  u(W_t, T-t)\cdot d(AW)_t \right)\times \right. \\
&   & \hspace*{12em}
      \times  \left. \left(\int_{0}^{T}
      \nabla  v(W_t, T-t)\cdot dW_t\right)\right]  \\
& = & (2\pi T)^{d/2}\expec 
      \int_{[0,T]} \langle 
      A^T\nabla  u(W_{t},T-t), \nabla  v(W_{t},T-t)
      \rangle dt \\
& = & \int_{\real^d\times [0,T]} \langle A^T\nabla  u(x,t),\nabla  v(x,t)
      \rangle (2\pi T)^{d/2} d \mu_{t,T}(x) dt,
\end{eqnarray*}
where $\mu_{t,T}= law (W_{T-t})$. Because of 
$\sup_{x\in\real^d, t>0} t^{\frac{d+1}{2}}|\nabla  u(x,t)| < \infty$
and similarly for $v$ and by splitting the integration over 
$\real^d\times [0,T]$ into $\real^d \times [0,T/2]$ and
$\real^d \times (T/2,T]$ we see by standard arguments that
\begin{multline*}
     \lim_{T\to\infty}
     (2\pi T)^{d/2}\expec \left[\left(\int_{0}^{T} 
     \nabla  u(W_t, T-t)\cdot d(AW)_t \right)\times \right. \\
     \times  \left. \left(\int_{0}^{T}
     \nabla  v(W_t, T-t)\cdot dW_t\right)\right]  
   = \int_{\real^d\times [0,\infty)} 
     \langle A^T\nabla u(x,t), \nabla  v(x,t) \rangle dx dt
\end{multline*}
(note that 
 $\int_{\real^d\times [0,\infty)} 
  |\nabla u(x,t)|^2 dx dt < \infty$ and the same for $v$, which 
follows by the argument below, cf. also 
\cite[Lemma 1.1] {Petermichel-Volberg}).
If $\F u$ is the Fourier transform of $u$ with respect to the 
$x$-variables, then it is well-known that 
$\F u (\xi ,t)= (\F f)(\xi )e^{-t|\xi |^2/2}$ for $\xi \in\real^d$ and 
$t>0.$ Observe also that $\F ((d/dx_k)u)=i\xi_k\F u$ for $k=1,...,d$.
By Parseval's formula and Fubini's theorem we may compute
\begin{eqnarray*}
&   & \int_{\real^d} (U_Af)(x)g(x)\, dx \\
& = & (2\pi)^{-d}\int_{\real^{d}}\left(\int_{\real^+}
      e^{-t|\xi |^2}\, dt\right) \F f(\xi )\overline{\F g(\xi )}
      \langle i\xi, A(i\xi)\rangle \, d\xi \\
& = & (2\pi)^{-d}\int_{\real^{d}} \F f(\xi )\overline{\F g(\xi )}
      |\xi|^{-2}\langle \xi,A\xi \rangle \,  d\xi \\
& = & \int_{\real^{d}} (T_mf)(x) g(x)\, dx ,
\end{eqnarray*}
where $m$ is the multiplier $m(\xi):= |\xi|^{-2}\langle \xi,A\xi \rangle$.
By  recalling that $R_j$ corresponds to
the multiplier $\xi_j/(i|\xi |)$ the claim follows immediately.
\halmos 
\end{proof}

%%%%%%%%%%%%%%%%%%%%%%%%%%%%%%%%%%%%%%%%%%%%%%%%%%%%%%%%%%%%%%%%%%%%%%%


\begin{thebibliography}{1}

\bibitem{Ban1}
R. Ba\~nuelos and G. Wang: {\it Sharp inequalities for martingales 
with applications to the Beurling-Ahlfors and Riesz transforms,} 
Duke Math. J.  80  (1995),  575--600. 

\bibitem{Ban2}
R. Ba\~nuelos and P. J. M\'endez-Hern\'andez: {\it Space-time Brownian 
motion and the Beurling-Ahlfors transform,} 
Indiana Univ. Math. J.  52  (2003),  no. 4, 981--990. 

\bibitem{Baernstein-MS}
Al. Baernstein and S. Montgomery-Smith:
{\it Some conjectures about integral means of $\partial f$ and $\partial^-f$,}
Complex Analysis and Differential Equations, edited by C.Kiselman, 
Acta Universitatis Upsaliensis C., Volume 64 (1999), 92-109.

\bibitem{Bass}
R.F. Bass:
{\em Probabilistic Techniques in Analysis,}
Springer, 1995.

\bibitem{Bou}
J. Bourgain:
{\it Some remarks on {B}anach spaces in which martingale difference
sequences are unconditional,}
Ark. Mat., 21:163--168, 1983.

\bibitem{Bou2}
J. Bourgain: {\it Vector-valued singular integrals and the $H\sp 1$-BMO duality.} In:  Probability theory and harmonic analysis (Cleveland, Ohio, 1983),  1--19, Monogr. Textbooks Pure Appl. Math., 98, Dekker, New York, 1986.

\bibitem{Burk1}
D. L. Burkholder : { \it A geometrical characterization of Banach spaces in which martingale 
difference sequences are unconditional,}  
Ann. Probab.  9  (1981),  997--1011. 

\bibitem{Burk2}
D. L. Burkholder : {\it A geometric condition that implies the existence of certain 
singular integrals of Banach-space-valued functions,}
in: Conference on harmonic analysis in honor of Antoni Zygmund, Vol. I, II 
(Chicago, Ill., 1981),  270--286, Wadsworth Math. Ser., Wadsworth, Belmont, CA, 1983.

\bibitem{Burk5} 
D.L. Burkholder : {\it Boundary value problems and sharp inequalities for 
                   martingale transforms,}
Ann. Prob. 12 (1984), 647--702. 

\bibitem{Burk4}
D.L. Burkholder : {\it Explorations in martingale theory and its applications,}
Ecole d'Et\'{e} de Probabilit\'{e}s de Saint-Flour, XIX--1989,
Lect. Notes Math.  1464, 1--66, 1992, Springer. 
    
\bibitem{Burk3}
D. L. Burkholder : { \it  Martingales and singular integrals in Banach spaces,} in: 
Handbook of the geometry of Banach spaces, Vol. I,  233--269, North-Holland, Amsterdam, 2001.

\bibitem{CoifmanWeiss} R. R. Coifman and G. Weiss: 
{\it Transference Methods in Analysis,}
C.B.M.S. Regional Conference Series in Math. No. 31, Amer. Math. Soc., Providence, RI, 1976.

\bibitem{DefantM-PHD}
M. Defant:
{\em Zur vektorwertigen Hilberttransformation,}
PhD thesis, Universit\"at Kiel, 1986.

\bibitem{Geiss8.5}
S. Geiss:
{\it A counterexample concerning the relation between decoupling constants
and UMD-constants,}
Trans. Amer. Soc., 351(4):1355--1375, 1999.

\bibitem{Grafakos} L. Grafakos: {\it  Classical and modern Fourier analysis,} 
Pearson 2004.

\bibitem{GunVar}
R. F. Gundy and N. Varopoulos: {\it Les transformations de Riesz et les intÃ©grales 
stochastiques,}
C. R. Acad. Sci. Paris S\'er. A-B  289  (1979), A13--A16.

\bibitem{Hyt}
T.P. Hyt\"onen: {\it Aspects of probabilistic Littlewood-Paley theory
in Banach spaces,} 
Banach spaces and their applications in analysis,  de Gruyter, Berlin,
343--355, 2007. 

\bibitem{James}
R.C. James: {\it Super-reflexive Banach spaces,}
Can. J. Math., 5 (1972), 896--904. 

\bibitem{deLeeuw} K. de Leeuw: 
{\it On $L_p$ multipliers,} Ann. of Math 81 (1965), 364--379.

\bibitem{Maurey} 
B. Maurey: 
{\it Syst\`{e}me de {H}aar,} Seminaire Maurey--Schwartz, Ecole Polytechnique, Paris, 1974--1975.

\bibitem{Petermichel-Volberg}
S. Petermichl and A. Volberg:
{\it Heating of the Ahlfors-Beurling operator: weakly quasiregular maps on
    the plane are quasiregular,}
Duke Math. J., 112(2):281--305, 2003.

\bibitem{Pietsch-Wenzel}
A. Pietsch and J. Wenzel:
{\it Orthonormal systems and {B}anach space geometry,}
Cambridge University Press, 1998.

\bibitem{Stein2} 
E. M. Stein: {\it Singular integrals and differentiability properties of functions,} 
Princeton University Press 1970.

\bibitem{SteinWeiss} E. M. Stein and G. Weiss: 
{\it  Introduction to Fourier analysis on Euclidean spaces,} 
Princeton University Press 1971.

\bibitem{Volberg} 
A. Volberg and F. Nazarov: {\it Heat extension of the Beurling operator and estimates for 
its norm,} (Russian. Russian summary), Algebra i Analiz 15 (2003),  142--158; 
translation in St. Petersburg Math. J. 15 (2004),  563--573.

\bibitem{Weis}
L. Weis: {\it Operator-valued Fourier multiplier theorems and maximal $L\sb p$-regularity,}  
Math. Ann.  319  (2001),  735--758.

\end{thebibliography}
\end{document}